# An Effective Estimate on Betti Numbers

Lei Ni[1]

**Abstract**
We provide an effective estimate on the Betti numbers of the loop space of a compact manifold which admits a finite Grauert tube. It implies the polynomial estimate in Chen (arXiv:2101.04368, 2021) after taking the radius of the tube to infinity.

## 1 Introduction

For a compact Riemannian manifold $(M, g)$, we say that it admits a Grauert tube of radius $R$ if there exists a *canonical* complex structure on the disc bundle of radius $R$. It is denoted by $T^R M$. The $R$ is called the radius. Such a complex structure which is called the adapted complex structure (cf. Theorem A of [13]) is uniquely determined by the Riemannian structure of $(M, g)$ and makes $M$ a totally real embedded submanifold. The concept arises in the study of the homogenous complex Monge–Ampère equation [2, 12] and the complexification of a real analytic manifold [5]. Please refer to [9] and reference therein for more detailed descriptions of its development. Its existence is not automatic. In fact it was proved by Lempert (cf. Theorem 1.5 of [8]) that the existence of a Grauert tube implies that $g$ must be a real analytic metric. Hence it does put some constraint on the metric $g$ besides requiring that $M$ is real analytic. Conversely for a compact real analytic manifold $M$ and a real analytic metric there always exists a real $R_0 > 0$ such that $T^{R_0} M$ is a Grauert tube. Namely there exists a *canonical adapted* complex structure on the disc bundle $T^{R_0} M$. The precise definition of the adapted complex structure requires introducing additional notions which we defer to the next section. The main result of this note is

[1] ✉ Lei Ni
  leni@ucsd.edu

1  Department of Mathematics, University of California, San Diego, La Jolla, CA 92093, USA

**Theorem 1.1** *Let $(M, g)$ be a simply connected compact Riemannian manifold which admits a Grauert tube of radius $R$. Given any coefficient field $\mathbb{F}$, and a positive integer $k$, then there exists a constant $C > 0$, which is independent of $\mathbb{F}$ and $k$, such that*

$$\sum_{i=0}^{k-1} \dim H_i(\Omega(M), \mathbb{F}) \leq \frac{\omega_{n-1}}{Vol(M)} \int_0^{Ck} \left( \frac{\sinh(\frac{\pi}{2R}\sigma)}{\frac{\pi}{2R}} \right)^{n-1} d\sigma. \quad (1.1)$$

*Here $\Omega(M)$ denotes the space of continuous loops, $\omega_k$ is the volume of the $k$-dimensional unit sphere $\mathbb{S}^k$.*

Taking $R \to \infty$, one can recover the polynomial (in terms of $k$) estimate of Chen [3] in the case when $M$ admits a global Grauert tube.

**Corollary 1.2** *Let $(M, g)$ be a compact Riemannian manifold which admits a global Grauert tube. Then*

$$\sum_{i=0}^{k-1} \dim H_i(\Omega(M), \mathbb{F}) \leq \frac{Vol(B^n(1))}{Vol(M)} (Ck)^n. \quad (1.2)$$

It follows that $M$ is rational elliptic, by the result of Félix and Halperin (cf. page 110 of [10]), and Proposition 5.6 of [10]. This provides a positive answer to a special case of a conjecture of Bott (which asserts that any simply connected close manifold with nonnegative sectional curvature must be rationally elliptic, in view of the result of Lempert–Szöke below, since the existence of a global Grauert tube implies the nonnegativity of the sectional curvature).

**Theorem 1.3** (Lempert–Szöke) *Let $(M, g)$ be a simply connected compact Riemannian manifold which admits a Grauert tube of radius $R$. Then the sectional curvature of $(M, g)$ is bounded from below by $-\frac{\pi^2}{4R^2}$. In particular, if $R = \infty$, $(M, g)$ has nonnegative sectional curvature.*

For a fixed point $p \in M$ let $D(p, R)$ denote the set of vectors in the tangent space $T_p M$ satisfying $|v| < R$. Namely $D(p, R) = \{v \in T_p M \mid |v| < R$. Let $n(p, R, x)$ be the number of the pre-images of $x \in M$ in the tangent space $T_p M$ under the exponential map at $p$ which are inside $D(p, R)$. Namely $n(p, R, x) = \#\{v \in T_p M \mid |v| < R, \exp_p(v) = x\}$. A result of Gromov, via the Morse theory on the energy functional defined on the space of pathes, asserts the following useful estimate on a compact simply connected manifold $M$ (cf. [10], Theorem 5.10, formula (5.3) and the estimate above it on page 124, and Remark 5.28) for a regular value $x$ of the exponential map $\exp_p : T_p M \to M$

$$\sum_{i=0}^{k-1} \dim H_i(\Omega(M), \mathbb{F}) \leq n(p, Ck, x).$$

Note that the right hand side is finite away from a measure zero set. Here $C$ is a constant independent of $\mathbb{F}$ and $k$. Hence

$$\sum_{i=0}^{k-1} \dim H_i(\Omega(M), \mathbb{F}) \leq \frac{1}{Vol(M)} \int_M n(p, Ck, x) \, d\mu(x). \tag{1.3}$$

On the other hand, the area formula for Lipschitz maps (cf. [4], Theorem 1 of Ch3.3) implies that

$$\int_M n(p, T, x) \, d\mu(x) = \int_0^T \int_{\mathbb{S}^{n-1}} \mathrm{Jac}(\exp_p)|_{(\sigma, \theta)} \, d\theta \, d\sigma \tag{1.4}$$

provided that $T \geq \mathrm{Diam}(M)$. Here $(\sigma, \theta)$ is the polar coordinate of $T_p M$, and $\mathbb{S}^{n-1}$ is the unit sphere in $T_p M$. It is well known that $\mathrm{Jac}(\exp_p)|_{(\sigma, \theta)}$ can be computed via the square root of the determinant of the gram matrix of the $n-1$ normal Jacobi fields $J_1(\sigma), \cdots, J_{n-1}(\sigma)$ along a geodesic $\gamma(\sigma)$ with $\dot{\gamma}(0) = \theta$, such that $J_i(0) = 0$, $\dot{J}_i(0) = e_i$, for $1 \leq i \leq n-1$. Here after adding $e_n = \dot{\gamma}(0) = \theta$ (assuming that $\gamma$ is parametrized by the arc-length), $\{e_i\}_{1 \leq i \leq n}$ forms an orthonormal frame of $T_p M$. We denote also their parallel transport as $\{e_i\}$. Combining them we have

$$\int_M n(p, T, x) \, d\mu(x) = \int_0^T \int_{\mathbb{S}^{n-1}} \sqrt{\det(\langle J_i(\sigma), J_j(\sigma) \rangle)_{1 \leq i,j \leq n-1}} \, d\theta \, d\sigma. \tag{1.5}$$

Note that if we let $J_n(\sigma) = \sigma e_n(\sigma)$, the right hand side remains the same if we replace the $(n-1) \times (n-1)$ matrix with the $n \times n$ matrix $(\langle J_i, J_j \rangle)_{1 \leq i,j \leq n}$ and then compensate with a factor of $\frac{1}{\sigma}$ since $\langle J_i, J_n \rangle = 0$ for $1 \leq j \leq n-1$ and $|J_n|^2 = \sigma^2$.

Note that the volume comparison fails beyond the first conjugate locus, hence could not be applied directly otherwise the conjecture of Bott would have been known many years ago. Here we show that it can be estimated under the additional structure of the Grauert tube via the analytic continuation. One would naturally conjecture that Theorem 1.1 (or the estimate on $\det(\langle J_i(\sigma), J_j(\sigma) \rangle)_{1 \leq i,j \leq n-1}$) holds for a closed manifold with a sectional curvature lower bound $-\frac{\pi^2}{4R^2}$, but without any assumption on the existence of the Grauert tube (in particular no real analyticity assumption on the metric).

## 2 Proof of the Theorem

The proof utilizes the framework developed, and results obtained, in the important paper of Lempert and Szöke [9]. Below we need some results from that paper after some basics on geodesic flows on which one can find more detailed coverage in some excellent books e.g. [7, 10, 11].

The idea is to use the existence of a holomorphically immersed strip $\mathcal{S}_R = \{z = \sigma + \sqrt{-1}\tau \in \mathbb{C} \mid 0 < \tau < R\}$, whose closure passes $p$, to estimate the right hand side of (1.5) via a Fatou Lemma (for positive harmonic functions) by adapting the considerations of [9] to prove the estimate (1.1). We include some introductory materials (mostly from [7, 9, 11]) for the convenience of the readers.

Let $\tau_M : TM \to M$ denote the bundle projection $TM \to M$. Recall that for any $u = (x^i, \xi^i) \in TM$, with $p = \tau_M(u)$ and $(x^1, \cdots, x^n)$ being a local chart near $p$, $u = \xi^i \frac{\partial}{\partial x^i}$, and the tangent $\eta = (x^i, \xi^i, X^i, \eta^i) \in T_u TM$ (for a curve $u(t) = (x^i(t), \xi^i(t))$, $X^i = \frac{d}{dt} x^i(t)$, $\eta^i = \frac{d}{dt} \xi^i(t)$), the tangent map $D\tau_M : TTM \to TM$ can be expressed as $D\tau_M|_u(\eta) = (x^i, X^i)$. We also abbreviate $D\tau_M$ as $(\tau_M)_*$. The kernel of $(\tau_M)_*$ is called the *vertical subspace* of $T_u TM$ (denoted by $V_u$). We may identify the vertical tangent at $u$, given locally as $(x^i, \xi^i, 0, \eta^i)$, with $(x^i, \eta^i)$, an element in the tangent space $T_p M$. Such identification map is denoted as $\iota_u$.

The *connection map* $K : TTM \to TM$ is defined as follows. For $u = (x^i, \xi^i)$, given any $X \in T_p M$ define the horizontal lift of $X = (x^i, X^i)$ (namely $X = X^i \frac{\partial}{\partial x^i}$) as $(x^i, \xi^i, X^i, -\Gamma^i_{jk} \xi^j X^k)$. One can check that the lift is so defined that it is $(x^i, \xi^i, X^i, \frac{d\xi^i}{dt}(0))$, the derivative of the vector $\xi(t)$ with $\xi(0) = u$ along a curve $c(t)$ with $c(0) = p$, $c'(0) = X$, [1] such that $\frac{D}{dt} \xi(t) = 0$. Here $\frac{D(\cdot)}{dt}$ denotes the covariant derivative with respect to the Levi-Civita connection. Namely it is the local derivative of the parallel vector extension along the direction given by $X$. Clearly the definition does not depend on the choice of $c(t)$. For $\eta = (x^i, \xi^i, X^i, \eta^i) \in T_u TM$ with $u = (x^i, \xi^i)$, $K(\eta)$ is defined as $(x^i, \eta^i + \Gamma^i_{jk} \xi^j X^k)$. Clearly if $\eta$ is a horizontal lift of $X$, $K(\eta) = 0$. One can also check readily that if $\xi(t) \in TM$ is a vector field along $\alpha(t)$ with $\xi(0) = \xi$ and $\xi'(0) = \eta$ and $\alpha(t) \in M$ with $\alpha(0) = p$, $\alpha'(0) = X$ (noting that $\alpha(t)$ can be taken to be $\tau_M(\xi(t))$), $K(\eta) = \frac{D\xi(t)}{dt}|_{t=0}$ (see also the proposition below). The kernel of $K$ at $u$ is called the *horizontal subspace* $H_u$. Clearly $V_u \cap H_u = \{0\}$ and $T_u TM = V_u \oplus H_u$ by the dimension counting. The following (cf. Proposition 4.1 of Ch II [11]) is well known.

**Proposition 2.1** *For $y \in T_p M$ and $X$ a smooth vector field (which is viewed as a map $M \to TM$), $K(DX(y)) = \nabla_y X$.*

From the above discussion it is clear that the map $j_u : T_u M \to T_p M \times T_p M$, defined by $j_u(\eta) = ((\tau_M)_*(\eta), K(\eta))$ for $\eta \in T_u TM$ is a linear isomorphism.

Let $\phi_\sigma : TM \to TM$ be the geodesic flow defined as $\phi_\sigma(v) = \dot{\gamma}(\sigma)$, with $\dot{\gamma}(0) = v$, $p = \pi(v)$, $\gamma(0) = p$. In terms of the local coordinate expression above, locally $v = (x_0^i, X^i)$ with $p = (x_0^1, \cdots, x_0^n)$. Then $\phi_\sigma(v) = (\gamma^i(\sigma), \dot{\gamma}^i(\sigma))$ with $\gamma^i$ being the $i$-th coordinate of $\gamma$ and $\gamma^i(0) = x_0^i$, $\dot{\gamma}^i(0) = X^i$. For the real line $\mathbb{R}$ one can identify $\mathbb{C}$ with $T\mathbb{R}$ via the identification of $\sigma + \sqrt{-1}\tau$ with $\tau \frac{\partial}{\partial \sigma}$. For any smooth map $f : M_1 \to M_2$ between two manifolds, let $f_* : TM_1 \to TM_2$ be the associated differential map which sends $(x, X(x)) \to (f(x), df|_x(X(x)))$. By the definition, $\gamma_*(\sigma_1 + \sqrt{-1}\tau_1) = (\gamma(\sigma_1), \frac{d}{ds}|_{s=0} \gamma(\sigma_1 + \tau_1 s)) = (\gamma(\sigma_1), \tau_1 \dot{\gamma}(\sigma_1))$. Abusing the notation we have that $\gamma_*(\sigma_1 + \sqrt{-1} \cdot 0) = \gamma(\sigma_1)$, namely we write $(\gamma(\sigma_1), 0)$ simply as $\gamma(\sigma_1)$. As $\gamma$ runs among all geodesics, $\gamma_* : T\mathbb{R} \setminus \mathbb{R} \to TM \setminus M$ defines a foliation, which is called the *Riemann foliation*. The uniqueness of the geodesic with a given initial point and a velocity vector implies that it is indeed a foliation.

---

[1] Here (˙) denotes the derivative with respect to the parameter for geodesics and ( )′ denotes the derivative for other parameters.

The adapted holomorphic structure can be defined as follows (cf. [13]).

**Definition 2.2** Let $(M, g)$ be a complete Riemannian manifold. For given $R$, a smooth complex structure on the manifold $T^R M$ will be called an adapted holomorphic structure if for any geodesic $\gamma : \mathbb{R} \to M$, the map $\gamma_* : \gamma_*^{-1}(T^R M) \to T^R M$ is holomorphic. When such adapted holomorphic structure exists the disc bundle $T^R M$ is called the Grauert tube of radius $R$.

Given the adopted holomorphic structure, an important object is the so-called *parallel vector field* (perhaps a more suitable name is the *generalized Jacobi field*) associated with $\gamma_*$. It is defined as the variational vector field of $(\gamma_t)_*$ for a family of geodesics $\gamma_t : \mathbb{R} \to M$, namely is defined as $\xi = \frac{D}{dt}|_{t=0}((\gamma_t)_*)$. Given a geodesic $\gamma(\sigma)$ with $\gamma'(0) = u$, and $\eta \in T_{\gamma_*(\tau_1 \sqrt{-1})} TM = T_{\tau_1 \dot\gamma(0)} TM$ for $\tau_1 \ne 0$, it is not hard to see that there exists a generalized Jacobi field (parallel vector field) $\xi$ such that $\xi(\sqrt{-1}\tau_1) = \eta$. (Since $\gamma$ is not parametrized by the arc length it is sufficient to consider $\tau_1 = 1$.) In fact if $\eta \in T_u TM$, let $z(t)$ be a path in $TM$ such that $z(0) = \dot\gamma(0) = u$, $z'(0) = \frac{\eta}{\tau_1}$. Now let $\gamma_t(\sigma)$ be the geodesic with initial condition $z(t)$, namely $\gamma_t(\sigma) = \tau_M(\phi_\sigma(z(t)))$. It is easy to see that $\xi(\sigma + \sqrt{-1}\tau) = \frac{D}{dt}|_{t=0}(\gamma_t)_*$ is a parallel field satisfying that $\xi(\sqrt{-1}\tau_1) = \eta$. In deed, $\frac{d}{dt}|_{t=0}(\tau_1 \phi_0(z(t))) = \tau_1 z'(0) = \eta$. The parallel vector field generalizes the concept of the Jacobi field, since by the way it is defined $\xi(\sigma) = \xi(\sigma + 0 \cdot \sqrt{-1}) = \frac{d}{dt}|_{t=0}(\gamma_t(\sigma), (0 \cdot \dot\gamma_t(\sigma)) = \frac{d}{dt}|_{t=0}\gamma_t(\sigma)$ (in the last equation we omit the second component since it is zero), which is the variational vector field of a family of geodesics, hence a Jacobi field. Namely $J(\sigma) = \xi|_\mathbb{R}(\sigma)$ is a Jacobi field. The Jacobi field is defined along the image of a geodesic $\gamma$, the generalized Jacobi field (parallel vector field) is defined along the image of $\gamma_*$, namely a leave of the Riemann foliation.

In the above if we choose $\tau_1 = 1$, $\tau_M(\phi_\sigma(z(t)))$ is a family of geodesic $\gamma_t(\sigma)$. Hence $J_u(0) = \frac{d}{dt}|_{t=0}\tau_M(z(t)) = d\tau_M(\eta)$. Since $\phi_0(z(t)) = z(t)$ is a path in $TM$, which covers $\alpha(t) = \tau_M(z(t))$, we have that $K(\eta) = \frac{Dz(t)}{dt}|_{t=0}$. On the other hand $\frac{DJ_u(\sigma)}{d\sigma}|_{\sigma=0} = \frac{D}{d\sigma}\frac{D}{dt}(\gamma_t(\sigma))|_{(0,0)} = \frac{D}{dt}\frac{D}{d\sigma}(\tau_M(\phi_\sigma(z(t)))|_{(0,0)} = \frac{D}{dt}|_{t=0}(\dot\gamma_t(0)) = \frac{Dz(t)}{dt}|_{t=0}$. Hence $K(\eta) = \frac{DJ_u}{d\sigma}(0)$ (cf. Lemma 4.3 of Ch 2 of [11]). The discussion works similarly for any $\tau_1 > 0$. The holomorphicity of $(\gamma_t)_*$ for a family of geodesics implies that the holomorphic component (namely the $(1, 0)$-part) of the variational vector field, the generalized Jacobi field (parallel vector field), is holomorphic by calculations with respect to the holomorphic coordinates (cf. Proposition 5.1 of [9]).

Now for an orthonormal frame $\{v_j\}_{j=1}^n$ as the above we choose $\xi_j$ and $\eta_j$ at $u \in TM$ such that $d\tau_M(\xi_j) = v_j$ and $K(\xi_j) = 0$ and $d\tau_M(\eta_j) = 0$, $K(\eta_j) = v_j$ and then extend them as above into $2n$ parallel vector fields along a Riemann foliation $\gamma_*(\sigma + \sqrt{-1}\tau)$. By Lemma 5.1 of [9], the holomorphic parts, $\xi_j^{1,0}$ and $\eta_j^{1,0}$ are holomorphic over the domain where $\gamma_*$ is holomorphic. In the case that $\gamma(\sigma)$ is parametrized by the arc-length, they are holomorphic on the strip $S_R = \{\sigma + \sqrt{-1}\tau \mid 0 < \tau < R\}$ if $T^R M$ is a Grauert tube with the adapted holomorphic structure. Here we may choose $\tau_1$ small such that $u \in T_p^R M$ in the construction given in the previous two paragraphs. The following proposition summarizes the main construction of [9].

**Proposition 2.3** *Let* $\Xi = (\xi_1, \cdots, \xi_n)$ *and* $H = (\eta_1, \cdots, \eta_n)$ *and* $\Xi^{1,0}$ *and* $H^{1,0}$ *be the holomorphic components. Then*

(i) *The 2n-vectors* $\{\xi_j, J(\xi_j)\}_{j=1}^n$ *are linearly independent on the strip* $\mathcal{S}_R$ *and* $\{\xi_j^{1,0}\}$ *are linearly independent over* $\mathbb{C}$ *on the strip (here J is the almost complex structure of the tube and* $\xi_j^{1,0} = \frac{1}{2}(\xi_j + \sqrt{-1}J\xi_j)$*); The 2n-vectors* $\eta_1, \cdots, \eta_n, \xi_1, \cdots, \xi_n$ *are linearly independent; Their restrictions to* $\mathbb{R}$ *are Jacobi fields* $J_1(\sigma), \cdots, J_n(\sigma), J_{n+1}(\sigma), \cdots, J_{2n}(\sigma)$ *satisfying that for* $1 \leq j \leq n$, $J_j(\sigma)|_{\sigma=0} = \eta_j|_\mathbb{R}(0) = 0, \frac{D}{d\sigma}J_j(\sigma)|_{\sigma=0} = \frac{D}{d\sigma}\eta_j|_\mathbb{R}(\sigma)|_{\sigma=0} = v_j$ *and* $J_{n+j}(\sigma)|_{\sigma=0} = \xi_j|_\mathbb{R}(0) = v_j, \frac{D}{d\sigma}J_{n+j}(\sigma)|_{\sigma=0} = \frac{D}{d\sigma}\xi_j|_\mathbb{R}(\sigma)|_{\sigma=0} = 0;$

(ii) *There exist a holomorphic matrix f on the strip* $\mathcal{S}_R$, *which, after being extended to* $\mathbb{R}$, *may have poles on a discrete subset of* $\mathbb{R}$, *such that* $H^{1,0}(z) = \Xi^{1,0}(z) f(z)$ *with* $z = \sigma + \sqrt{-1}\tau \in \mathcal{S}_R$;

(iii) Im(*f*) *is symmetric and positive definite for* $\sigma + \sqrt{-1}\tau$ *with* $\tau > 0$; *On* $\mathbb{R}$, *f is real and* $H = \Xi f$ *over the set where f is finite.*

(iv) $f|_\mathbb{R}$ *is symmetric over where it is finite, f is holomorphic near 0, and* $f(0) = 0$ *and* $f'(0) = \mathrm{id}$.

*Proof* The results were mainly proved in Sect. 6 of [9] by exploiting the Kähler/symplectic structure on $T^R M$. By the way in which $\xi_j$ and $\eta_j$ are constructed it is clear that they are linearly independent at $u \in TM$ due to that $j_u$ is an isomorphism. Then their extensions are linearly independent (over $\mathbb{R}$) due to the property of the geodesic flow $\phi_\sigma$, in particular, it is an isomorphism between tangent spaces of the domain and target points (cf. Propositions 1.90 and 1.92 of [1]). The linear independence of $\{\xi^{1,0}\}$ over $\mathbb{C}$ needs to use the Kähler form on $T^R M$, which was proved in Proposition 6.4 of [9]. The positivity of Im $f$ is proved in Lemma 6.7 of [9]. The holomorphicity of $f$ near 0 is due to the fact that $\{\xi_j(\sigma)\}$ are linearly independent for $\sigma$ small. The equations satisfied by $f(0)$ and $f'(0)$ are easy consequences of the constructions of $\{\xi_j\}$ and $\{\eta_j\}$. □

Applying the arguments/proofs in [9] we also have the following result.

**Proposition 2.4** *Under the assumption that M admits a Grauert tube of radius R we have the following results.*

(i) *The 2n-vectors* $\{\eta_j, J(\eta_j)\}_{j=1}^n$ *are linearly independent on the strip* $\mathcal{S}_R$ *and* $\{\eta_j^{1,0}\}$ *are linearly independent over* $\mathbb{C}$ *on the strip;*

(ii) *There exist a holomorphic matrix* $\tilde{f}$ *on the strip* $\mathcal{S}_R$, *which, after being extended to* $\mathbb{R}$, *may have poles on a discrete subset of* $\mathbb{R}$, *such that* $\Xi^{1,0}(z) = H^{1,0}(z) \tilde{f}(z)$ *with* $z = \sigma + \sqrt{-1}\tau \in \mathcal{S}_R$; $\tilde{f} = f^{-1}$, *hence symmetric over the points where* $\tilde{f}$ *is finite;*

(iii) Im($\tilde{f}$) *is symmetric and negative definite for* $\sigma + \sqrt{-1}\tau$ *with* $\tau > 0$; *On* $\mathbb{R}$, $\tilde{f}$ *is real and* $\Xi = H \tilde{f}$, *provided that* $\tilde{f}$ *is finite; And* $\tilde{f}$ *has a pole at* $z = 0$ *and is finite for* $\sigma \neq 0$ *small. Moreover* $\tilde{f}$ *extends to* $\mathbb{R} \setminus S$ *with S being the set* $\{\sigma_j\} \cup \{0\}$ *where* $\{\gamma(\sigma_j)\}$ *is the set of the conjugate points with respect to* $\gamma(0) = p$.

*Proof* First by Proposition 6.6 and the proof of Proposition 6.4 in [9] we have that $\{\eta_j^{1,0}(z)\}$ are linearly independent for $\tau > 0$. Hence there exists a matrix valued holomorphic function $\tilde{f}$ such that $\Xi^{1,0}(z) = \mathrm{H}^{1,0}(z)\tilde{f}(z)$. Namely $f(z)$ is invertible with $\tilde{f} = f^{-1}$ for $\tau > 0$. Hence $\tilde{f}$ and $\mathrm{Im}\,\tilde{f}$ are symmetric. This proves (i) and (ii).

By the exactly same argument as that of Proposition 6.8 of [9], $\mathrm{Im}\,\tilde{f}$ is invertible for $\tau > 0$. Now we consider the expansion of $\tilde{f}$ near 0:

$$\tilde{f}(\sqrt{-1}\tau) = \left(f(0) + f'(0)(\sqrt{-1}\tau) + O(\tau^2)\right)^{-1} = \left((\sqrt{-1}\tau(\mathrm{id} + O(\tau))\right)^{-1}$$
$$= \frac{-\sqrt{-1}}{\tau}\mathrm{id} + O(1).$$

This implies that $\mathrm{Im}(-\tilde{f})$ is positive definite at $z = \sqrt{-1}\tau$ with $\tau$ small. The analyticity of $\tilde{f}$ and the fact that $\mathrm{Im}\,\tilde{f}$ is invertible for $\tau > 0$ imply that $\mathrm{Im}(-\tilde{f})$ is positive definite for $z \in \mathcal{S}_R$. This proves the first part of (iii).

Since on $\mathbb{R}$, $f$ is real and $\Xi = \mathrm{H}\,\tilde{f}$ over the points where $\tilde{f}$ is finite. Hence $\tilde{f}$ does not have a pole, except possibly at $\sigma_j$, where $\gamma(\sigma_j)$ is conjugate to $p$ along $\gamma$ from the construction of $\xi_j(\sigma)$ and $\eta_j(\sigma)$. Since $\{\eta_j|_{\mathbb{R}}(\sigma)\}$ are zero at $\sigma = 0$, and are linearly independent for $\sigma$ small, and being linearly independent for $\sigma$ if $\gamma(\sigma)$ is not conjugate to $p$ along $\gamma$, the last part of (iii) holds. □

Note that both $f(z)$ and $-\tilde{f}(z)$ are valued in the Siegel upper-half space of degree $n$. The following result provides the key ingredient for estimating the Jacobian of the exponential map.

**Proposition 2.5** *For all $\sigma \in \mathbb{R}$, except a discrete subset,*

$$\Xi^{\mathrm{tr}}\mathrm{H}' - (\Xi^{\mathrm{tr}})'\mathrm{H} = \mathrm{id}; \quad -(\tilde{f}^{\mathrm{tr}})'\mathrm{H}^{\mathrm{tr}}\mathrm{H} = \mathrm{id}, \quad \Xi^{\mathrm{tr}}\Xi f' = \mathrm{id}. \quad (2.1)$$

*Here we identify the $\Xi(\sigma)$ and $\mathrm{H}(\sigma)$ with their matrices representation, with respect to a orthonormal frame $\{e_i\}$ obtained by parallel transplanting $\{e_i\}$, a frame at $p$, along $\gamma(\sigma)$.*

**Proof** First observe that for a holomorphic $F(z)$ with $z = \sigma + \sqrt{-1}\tau$ defined on $\mathcal{S}_R$, if $F = U + \sqrt{-1}V$, with $U$ being the real part and $V$ being the imaginary part, the Cauchy–Riemann equation implies that $F_z = U_\sigma + \sqrt{-1}V_\sigma$. In the case that $F|_{\mathbb{R}}$ is real we have that $F_z = U_\sigma = F_\sigma$. Namely by abusing the notation we denote by $F'$ both the complex derivative and the $\frac{d}{d\sigma}$ when $F$ has a finite extension near a point $\sigma_0 \in \mathbb{R}$. With the above consideration $\Xi_\sigma = 2\Re(\Xi_z^{1,0}) = 2\Re(\mathrm{H}_z^{1,0}\tilde{f} + \mathrm{H}^{1,0}\tilde{f}_z)$ on $\mathcal{S}_R$ and over the domain where $\tilde{f}$ has a holomorphic extension. Restrict to $\mathbb{R} \setminus S$, where $\tilde{f}$ is finite, we have that the right hand side equals to $\mathrm{H}_\sigma\tilde{f} + \mathrm{H}\tilde{f}_\sigma$ since $\tilde{f}$ is real valued there. Namely $\Xi_\sigma = \mathrm{H}_\sigma\tilde{f} + \mathrm{H}\tilde{f}_\sigma$ holds on $\mathbb{R} \setminus S$, which we abbreviate as $\Xi' = \mathrm{H}'\tilde{f} + \mathrm{H}\tilde{f}'$. The second identity follows from the first by plugging $(\Xi^{\mathrm{tr}})' = (\tilde{f}^{\mathrm{tr}})'\mathrm{H}^{\mathrm{tr}} + \tilde{f}^{\mathrm{tr}}(\mathrm{H}^{\mathrm{tr}})'$ into the first equation, and noting that $\mathrm{H}^{\mathrm{tr}}\mathrm{H}' - (\mathrm{H}^{\mathrm{tr}})'\mathrm{H} = 0$ by Proposition 6.10 of [9] (the argument below also provides a simple proof).

The first identity follows from that (i) $\frac{d}{d\sigma}(\Xi^{\mathrm{tr}}\mathrm{H}' - (\Xi^{\mathrm{tr}})'\mathrm{H}) = 0$, which is a consequence of the Jacobi equation; and (ii) $(\Xi^{\mathrm{tr}}\mathrm{H}' - (\Xi^{\mathrm{tr}})'\mathrm{H})(0) = \mathrm{id}$. This argument also provides a simple proof of Proposition 6.10 of [9].

The last identity was proved in the proof of Proposition 6.11 of [9] for $\sigma$ small. To see that it holds on $\mathbb{R} \setminus S'$ with $S'$ being the set of the poles of $f$ we can apply the same argument as above. Namely plug $(\mathrm{H})' = \Xi' f + \Xi f'$, which holds on $\mathbb{R} \setminus S'$, into the first identity, and note that $\Xi^{\mathrm{tr}} \Xi' - (\Xi^{\mathrm{tr}})' \Xi = 0$ on $\mathbb{R}$, we have the last identity for all $\sigma \in \mathbb{R} \setminus S'$. □

By a calculation similar to the one in the proof of Proposition 7.1 of [9], the Jacobi curvature $R_\gamma(v) = R(v, \dot{\gamma})\dot{\gamma}$ has the expression in terms of the Schwarz derivative of $\tilde{f}$ over the $\sigma \in \mathbb{R} \setminus S$, with $S$ being the set of the pole of $\tilde{f}$ (which contains $\{\sigma_j\}$ with $\{\gamma(\sigma_j)\}$ being the conjugate points of $\gamma(0) = p$ along $\gamma$):

$$R_\gamma = \frac{1}{2} \mathrm{H} \left( \tilde{f}'''(\tilde{f}')^{-1} - \frac{3}{2}(\tilde{f}''(\tilde{f}')^{-1})^2 \right) \mathrm{H}^{-1}. \tag{2.2}$$

The result can also be derived from the invariance of the Schwarz derivative under the projective transformation and Proposition 7.1 of [9]. Namely if we denote the Schwarz derivative of $f$ as $\{f\}$, $\{\tilde{f}\} = \tilde{f}\{f\}\tilde{f}^{-1}$. Recall that for any nonsingular family of matrices $A(\sigma)$, the Cartan matrix $C(A)$ is defined as $A'A^{-1}$. The Schwarz derivative $\{f\}$ can be expressed as $C^2(f') - \frac{1}{2}(C(f'))^2$.

**Corollary 2.6** *Away from a discrete subset $S \subset \mathbb{R}$,*

$$(\langle J_j(\sigma), J_k(\sigma) \rangle) = f(\sigma)(f'(\sigma))^{-1} f(\sigma).$$

*In particular*

$$\det(\langle J_j(\sigma), J_k(\sigma) \rangle) = \frac{1}{\det(-\frac{d}{d\sigma}(f^{-1}(\sigma)))}. \tag{2.3}$$

*The conjugate points of $p = \gamma(0)$ are at $\gamma(\sigma_j)$ with $\sigma_j$ being the poles of $\det(\tilde{f}')$.*

*Proof* The result, namely the first equation in the corollary, follows from Proposition 2.5, since it implies that $\mathrm{H}^{\mathrm{tr}} \mathrm{H} = \left(-\tilde{f}'\right)^{-1}$ and $(\langle J_j(\sigma), J_k(\sigma) \rangle)$ is nothing but the matrix representation of $\mathrm{H}^{\mathrm{tr}} \mathrm{H}$. The set $S$ can be chosen to be the union of poles of $\tilde{f}$ and zeros of $\det(-(\tilde{f}^{\mathrm{tr}})')$. The Eq. (2.3) follows by taking the determinant on both sides.

The main equation also follows from a less direct argument (cf. [3]) as follows. If $\gamma(\sigma_0)$ is not a conjugate point, $\mathrm{H}(\sigma_0)$ is non-degenerate. Hence $\tilde{f}(\sigma_0) = f^{-1}(\sigma_0)$ is finite. The Morse index theorem implies that away from a discrete $S$, $\tilde{f}(\sigma) = f^{-1}(\sigma)$ is analytic. Let $\{e_j\}$ be the standard orthonormal basis of $\mathbb{R}^n$. Then

$$\begin{aligned}
\langle J_j(\sigma), J_k(\sigma) \rangle &= \langle \mathrm{H}(\sigma) e_j, \mathrm{H}(\sigma) e_k \rangle \\
&= \langle \Xi(\sigma) f(\sigma) e_j, \Xi(\sigma) f(\sigma) e_k \rangle \\
&= \langle \Xi(\sigma) f(\sigma) e_j, \Xi(\sigma) f(\sigma) e_k \rangle \\
&= \langle f(\sigma)(\Xi(\sigma))^{\mathrm{tr}} \Xi(\sigma) f(\sigma) e_j, e_k \rangle.
\end{aligned}$$

This shows that $(\langle J_j, J_k \rangle) = f(\sigma)(\Xi(\sigma))^{\text{tr}}\Xi(\sigma)f(\sigma)$. By Proposition 2.5, $(\Xi(\sigma))^{\text{tr}}\Xi(\sigma) = (f'(\sigma))^{-1}$ for $\sigma \in \mathbb{R} \setminus S$. Putting them together we have the claim. □

This reduces the estimate of the right hand side of (1.5) to the estimate of the integral of the quantity in the right side of (2.3), which has a holomorphic extension on the strip $\mathcal{S}_R$. We also note that $\tilde{f}'$ is always invertible on $\mathbb{R} \setminus S$, given that H has a finite extension on $\mathbb{R}$.

The following Fatou type lemma plays a crucial role in both [9] (in proving Theorem 1.3) and [3], whose proof is reduced to a representation formula for positive harmonic functions on $\mathbb{C}^+$ (cf. Theorem 3.1.8 of [6]).

**Proposition 2.7** *Let $F(\zeta)$ be a $n \times n$ matrix valued holomorphic function on the upper half plane $\mathbb{C}^+ = \{\zeta \mid \text{Im}\,\zeta > 0\}$ that has a holomorphic extension $\mathbb{R}$ except a discrete set which are the poles of the analytic continuation of $F$. Then there exists a $n \times n$ symmetric matrix of Borel measures $\mu = (\mu_{jk})$ defined on $\mathbb{R}$ such that*

(i) $\operatorname{supp}\mu_{jk}$ *contains only a discrete set;*
(ii) $\int_{-\infty}^{\infty} \frac{|d\mu_{jk}(t)|}{1+t^2} < \infty$;
(iii) $\mu$ *is positive definite in the sense that $\langle \mu(w), w \rangle$ is a positive measure if $w \in \mathbb{R}^n$;*
(iv) $\frac{\partial}{\partial \zeta}(F(\zeta)) = A + \frac{1}{\pi}\int_{-\infty}^{\infty} \frac{d\mu(t)}{(\zeta-t)^2}$, *for $\zeta \in \mathbb{C}^+$. Here A is symmetric positive semidefinite matrix.*

**Proof** Please see Proposition 7.4 of [9]. The trick is to define $u_\xi = \text{Im}\langle F\xi, \xi \rangle$ for any $\xi \in \mathbb{R}^n$ and apply the result (Riesz-Herglotz representation theorem, cf. Theorem 3.1.8 of [6]) for nonnegative harmonic functions. This gives a measure $\mu_\xi$ which is a quadratic form of $\xi$, hence a symmetric matrix of measures which is nonnegative in the sense of (iii). Part (ii) is obtained by the conformal transformation from the unit disc to the upper half plane. The part (i) is due to that the measure $\mu_\xi$ is the weak limit of $u_\xi$. One obtains an expression for $\frac{\partial}{\partial \zeta}(F(\zeta))$ by expressing the derivative in terms of the derivative of imaginary part of $\tilde{F}$ (as in the proof of Proposition 2.5). Part (iv) also follows from [6], Exercise 3.1.6, which asserts that

$$F(\zeta) = A\zeta + B + \frac{1}{\pi}\int_{-\infty}^{\infty}\left(\frac{-1}{\zeta - t} - \frac{t}{1+t^2}\right)d\mu(t),$$

by taking derivative on both sides. Here $B$ is a constant matrix. □

Now let $F(\zeta) = -\tilde{f}(\frac{R}{\pi}\log\zeta) = -f^{-1}(\frac{R}{\pi}\log\zeta)$. For $\zeta \in \mathbb{C}^+$, $z = \frac{R}{\pi}\log\zeta$ is inside the strip $\mathcal{S}_R$. Hence $F(\zeta)$ is defined on $\mathbb{C}^+$. By Proposition 2.4, $F(\zeta)$ satisfies the assumptions of Proposition 2.7.

Applying Proposition 2.7 to $F(\zeta)$ we have the associated matrix of Borel measures $(\mu_{jk})$. Note that $\tilde{f}(\sigma)$ only have discrete poles on $\mathbb{R}$ at $\{\sigma_1, \sigma_2, \cdots\}$. As $f(0) = 0$, 0 is a pole of $\tilde{f}(\sigma)$. Hence $F(\zeta)$ has a pole at $\zeta = 1$. Since $\tilde{f}(\sigma)$ has only discrete poles, $\mu$ has support at discrete points $\{t_1, t_2, \cdots\}$

$$F'(t) = A + \frac{1}{\pi}\sum_j \frac{\mu(t_j)}{(t-t_j)^2}, \text{ for } t \in \mathbb{R}\setminus\{t_1, t_2, \cdots\}. \quad (2.4)$$

Here $A$ is a semidefinite positive constant matrix. Note $F'(t) = (-f^{-1})' \left( \frac{R}{\pi} \log t \right) \frac{R}{\pi t}$ with $\sigma = \frac{R}{\pi} \log t$. We shall show that

$$\frac{1}{\det(-\frac{d}{d\sigma}(\tilde{f}(\sigma)))} \le \left( \frac{R^2}{\pi^2} \right)^n \left( e^{\frac{\sigma\pi}{R}} + e^{-\frac{\sigma\pi}{R}} - 2 \right)^n. \tag{2.5}$$

Theorem 1.1 then follows by a simple calculation.

Since the measure $\mu$ is the weak limit of Im $F\, dt$ with $\zeta = t + \sqrt{-1}s$ (cf. Theorem 3.1.8 of [6]), calculations similar to [3] show that

$$\begin{aligned}
\mu(1) &= \lim_{\delta \to 0} \mu((1-\delta, 1+\delta)) = \lim_{\delta \to 0} \lim_{s \to 0} \int_{1-\delta}^{1+\delta} \mathrm{Im}(-f^{-1}(\frac{R}{\pi} \log(t + \sqrt{-1}s)))\, dt \\
&= \lim_{\delta \to 0} \lim_{s \to 0} \int_{1-\delta}^{1+\delta} -\mathrm{Im} \left( f(0) + f'(0)\frac{R}{\pi}(t + \sqrt{-1}s - 1) + O(|t + \sqrt{-1}s - 1|^2) \right)^{-1} dt \\
&= \frac{\pi}{R} \lim_{\delta \to 0} \lim_{s \to 0} \int_{1-\delta}^{1+\delta} -\mathrm{Im} \left( \frac{1}{t + \sqrt{-1}s - 1} \mathrm{id} \right) + O(1)\, dt \\
&= \frac{\pi}{R} \lim_{\delta \to 0} \lim_{s \to 0} \int_{1-\delta}^{1+\delta} \frac{s}{(t-1)^2 + s^2} \mathrm{id}\, dt = \frac{\pi^2}{R} \mathrm{id}.
\end{aligned}$$

This, together with (2.4) and the positivity of the measure $\mu$, implies that

$$(-\tilde{f})'(\sigma) \frac{R}{\pi t} \ge \frac{\pi}{R} \frac{\mathrm{id}}{(t-1)^2}.$$

Taking determinant on the both sides and noting $t = e^{\frac{\pi\sigma}{R}}$, (2.5) follows from Proposition 2.5 and calculations.

Finally to get the estimate in Theorem 1.1 one simply observes that when restricted to $\mathbb{R}$, $f_{jn} = 0$, for $1 \le j < n$, and $f_{nn} = \sigma$, and applies the above argument to the up-left $(n-1) \times (n-1)$ sub-matrix of $f$.

**Acknowledgements** We would like to thank Burkhard Wilking for bringing my attention to the subject, László Lempert for helpful discussions on the subject and results of [8, 9].